\documentclass[11pt]{amsart}
\usepackage{amsmath,amsfonts,amssymb,amsthm}
\usepackage{graphics,subfig}
\usepackage{epsfig}
\usepackage{esint}
\usepackage{footmisc}
\usepackage{tikz}
\usepackage{enumerate}
\usepackage{aligned-overset}

\numberwithin{equation}{section}

	\DeclareMathOperator	{\E}			{\mathbb{E}} 
	\DeclareMathOperator	{\N}			{\mathbb{N}}
	\renewcommand	{\P}			{\mathbb{P}}
	
	\DeclareMathOperator	{\R}			{\mathbb{R}}
	\DeclareMathOperator	{\Z}			{\mathbb{Z}}
	
	\DeclareMathOperator	{\supp}			{supp}

\newcommand{\Wb}{\mathsf{W}\mathsf{b}}
\newcommand{\Wa}{\mathsf{W}\mathsf{a}}
\newcommand{\W}{\mathsf{W}}
\newcommand{\q}{ \mathsf{q}}

\newcommand{\ignore}[1]{}

\usepackage[hidelinks,hyperindex,breaklinks]{hyperref}

\newtheorem{theorem}{Theorem}[section]

\newtheorem{proposition}[theorem]{Proposition}

\newtheorem{remark}[theorem]{Remark}
\newtheorem{lemma}[theorem]{Lemma}
\newtheorem{corollary}[theorem]{Corollary}

\usepackage[margin=1.5in]{geometry}

\title[Asymptotic transportation cost]{Asymptotic transportation cost of the optimal matching subject to different boundary conditions}

\author{Martin Huesmann, Arthur Wachtel}

\begin{document}
\begin{abstract}
    We establish equality between several asymptotic limits of the transportation cost between stationary random measures, with respect to the Euclidean distance cost raised to the $p$-th power, with different boundary conditions. 

    Our approach relies on a comparison with the optimizer of the mean transportation cost between the two random measures, which is defined on the whole space. In particular we do not use any PDE techniques related to the regularity theory for optimal transport maps.
\end{abstract}
\maketitle
\section{Introduction}
The goal of this note is to answer a question from \cite{GHT} regarding the equality of some constants appearing in thermodynamic limits of the random optimal matching problem with different boundary conditions. This is an instance of random combinatorial optimization problems which has been at the center of attention for a very long time due to its connections to computer science and statistical physics.
\medskip

It is most easily stated in the following way: Given two independent point clouds $(X_i)_{i=1}^n$ and $(Y_i)_{i=1}^n$ which are distributed according to some distribution $\lambda$ on $\R^d$ and a parameter $p\ge 1$, we are interested in 
\begin{align}\label{matching}
    \inf_{\pi\in\mathcal{S}_n}\frac 1n\sum_{i=1}^n|X_i-Y_{\pi(i)}|^p,
\end{align}
where $\mathcal{S}_n$ is the set of all permutations of $\{1,...,n\}$. By the Birkhoff -- von Neumann Theorem, we might reformulate this in terms of the Wasserstein $p$-cost $\W_p^p$. Indeed, writing $\mu=\frac1n\sum_{i=1}^n\delta_{X_i}$ and $\nu=\frac1n\sum_{i=1}^n\delta_{Y_i}$,
\begin{align*}
     \inf_{\pi\in\mathcal{S}_n}\frac 1n\sum_{i=1}^n|X_i-Y_{\pi(i)}|^p=\W^p_p(\mu,\nu).
\end{align*}
This identity allows a more general framework, for example point clouds with different cardinality or non-atomic measures. Especially the case where $\mu$ is kept as before, while $\nu$ is replaced by the law $\lambda$, has been investigated a lot in recent history. On the one hand it is particularly amenable to treatment using the regularity theory for optimal transport maps (see e.g. \cite{AGlT, CM, GHO, Koc}), while also being highly relevant in applications, since it can be interpreted as a quantization cost \cite{DSS,DY}. The main result of this viewpoint is the following high resolution formula, which appears in this form only in \cite{GHT}, even though different, weaker formulations can be found in \cite{BB,DSS,DY}.

\begin{theorem}\label{T:DSS}
    There exist constants $\underline{c}(p,d)\le \overline{c}(p,d)$ such that for any $d\ge 3, p\ge1$ and every probability measure $\lambda$ on $\R^d$, assuming either
    \begin{itemize}
        \item $p<\frac d2$ and $\lambda$ has finite $q$-th moment for some $q>\frac {dp}{d-p}$
        \item $\lambda$ is supported on a smooth, bounded, and connected set and is absolutely continuous with respect to the Lebesgue measure $\mathcal{L}$ and has H\"older continuous density bounded from above and below on its support
    \end{itemize}
    it holds that
    \begin{align}\label{eq:highresolution}
        \underline{c}(p,d)\int_{\R^d} \lambda_a^{1-p/d}&\le \liminf_{n\to\infty}n^{p/d}\E[\W^p_p(\frac1n\sum_{i=1}^n\delta_{X_i},\lambda)]\\
        &\le \limsup_{n\to\infty}n^{p/d}\E[\W^p_p(\frac1n\sum_{i=1}^n\delta_{X_i},\lambda)]\le \overline{c}(p,d)\int_{\R^d} \lambda_a^{1-p/d}. \nonumber
    \end{align}
    Here $\lambda_a$ denotes the absolutely continuous part of the measure $\lambda$ and $(X_i)_i$ is a sequence of i.i.d. random variables with law $\lambda$. The constants $\underline{c}(p,d),\overline{c}(p,d)$ are given by
    \begin{align*}
        \underline{c}(p,d)&:=\lim_{n\to\infty}n^{-d}\E[\Wb^p_{\Lambda_n}(\eta,\kappa_n\mathcal{L})]\\
         \overline{c}(p,d)&:=\lim_{n\to\infty}n^{-d}\E[\W^p_p(\eta\llcorner\Lambda_n,\kappa_n\mathcal{L}\llcorner\Lambda_n)],
    \end{align*}
    where $\Lambda_n=(-\frac{n}{2},\frac{n}{2})^d,$  $\eta$ is a homogeneous Poisson point process on $\R^d$, $\kappa_n:=\frac{\eta(\Lambda_n)}{n^d}$ and $\Wb^p_{\Lambda_n}$ denotes the boundary $p$-Wasserstein distance\footnote{For precise definitions we refer to Subsection \ref{SS:1} and \ref{SS:2}} on $\Lambda_n$.
\end{theorem}

The existence of the constants $\overline c(p,d)$ and $\underline c(p,d)$ can be deduced from now-standard sub- and super-additivity arguments which are typical in the study of Euclidean combinatorial optimization problems. For a proof specific to the optimal matching see \cite{GTa} and for a more general approach see \cite{BB, GTb}.
\medskip

A natural conjecture is of course $\underline{c}(p,d)=\overline{c}(p,d)$, which would also imply existence of the limit in Theorem \ref{T:DSS}. This question was partially resolved in \cite{GHT}, where a PDE approach similar to \cite{AGT} was employed to establish
\begin{align}\label{eq:GHT}
    \underline{c}(2,d)=\overline{c}(2,d)&&\text{for any $d\ge 3$.}
\end{align}

In \cite{DSS}, the lower bound in \eqref{eq:highresolution} is given by a smaller constant $c_a(p,d)$, which is defined through an asymmetric transport distance,
\begin{align*}
    c_a(p,d):=\lim_{n\to\infty}n^{-d}\E[\Wa_{\Lambda_n}^p(\eta,\kappa_n\mathcal{L})].
\end{align*}
The authors of \cite{GHT} raise the question whether $c_a(2,d)=\underline{c}(2,d)$, and thus for $p=2$, these two constants would also agree with $\overline{c}(2,d)$. We answer this question affirmatively, for any $p\geq 1$,  by relating both constants to another, namely the mean transportation cost $c_\infty(p,d)$, which was introduced in \cite{HS,Hue} and is defined by a Wasserstein distance on the space of laws of equivariant random measures (cf.\ \eqref{D:Global}). 
\begin{theorem}\label{T:2}
    For all $d\ge 3$ and $p\ge 1$,
    \begin{align*}
        \underline c(p,d)=c_\infty(p,d)=c_a(p,d).
    \end{align*}
\end{theorem}

Let us briefly comment on the dependence on the dimension. The randomness introduces local fluctuations, which leads to a  logarithmic divergence in dimension $d=2$, as was first realized in \cite{AKT}. This was subsequently used in \cite{HMOa, HMOb} to establish that no stationary $p$-cyclically monotone matching between a Poisson point process and the Lebesgue measure exists for any $p\in(1,\infty)$  and therefore $c_\infty(p,2)=\infty$. 
\medskip

As noted in \cite{GHT}, establishing Theorem \ref{T:2} suffices to deduce that all three transport distances in the limit $n\to\infty$
generate the coupling which optimizes the mean transportation cost. 
\begin{corollary}\label{C:1}
    Let $\{\q_n\}_n$ be a sequence of optimizers to any of the transport distances $\W_2^2(\eta\llcorner\Lambda_n,\kappa_n\mathcal{L}\llcorner\Lambda_n),\Wa_{\Lambda_n}^2(\eta,\kappa_n\mathcal{L})$, or $\Wb_{\Lambda_n}^2(\eta,\kappa_n\mathcal{L})$. Then stationarised versions of $\{\q_n\}_n$ converge vaguely to the unique invariant optimizer $\q$ which achieves the mean transportation cost $c_\infty(2,d)$.
\end{corollary}

This can be obtained by repeating the proof of \cite[Thm. 4.8]{HS}. The interesting phenomenon is that even though the boundary and the asymmetric variants are allowed to not use all the mass of $\eta$ within $\Lambda_n$, this does not leave any trace in the limit and in fact all the mass of $\eta$ is used.
\medskip

\begin{remark}
Theorem \ref{T:2} and Corollary \ref{C:1} are consequences of more general results where $\eta$ and $\mathcal L$ are replaced by jointly invariant random measures of equal intensity, see Theorem \ref{T:1}. In fact our approach uses neither any special properties of the Poisson process nor the fact that 
we match to a smooth, deterministic target measure.
\end{remark}

In the generality of the current article the four constants cannot agree, as they do not, for instance, in the case of a stationarised lattice and the Lebesgue measure. For sufficiently random random measures we do expect these four constants to agree for any $p,$ but this remains open for now.

\subsection{Related literature}
The boundary transport was first introduced in \cite{FG} to study parabolic equations subject to Dirichlet boundary conditions. However it has proven to be a very useful lower bound to the usual $p$-Wasserstein distance and has been employed to determine the limit of \eqref{matching} in the critical dimension $d=2$ when the unit cube is replaced by an arbitrary domain with Lipschitz boundary in \cite{AGT}, extending the seminal work \cite{AST}, which rigorously proved a challenging PDE ansatz from \cite{CLPS}. See also \cite{ATV} and \cite{AGGT} for an extension to arbitrary $1\le p<\infty$.

A byproduct of the recent result \cite{Xu} is the extension of \cite{GHT} to the case of matching of two independent Poisson processes in dimensions $d\geq 5,$ i.e.\ for point to point variants of $\underline{c}(2,d)$ and $\overline{c}(2,d)$ the equality corresponding to \eqref{eq:GHT} holds for $d\geq 5.$ Consequently, also in that situation all four constants are equal.

\section{Preliminaries}\label{S:2}
In this section we will define the various Wasserstein distances mentioned in the introduction as well as review some of their properties. Moreover we introduce the theory of optimal transport for (infinite) random measures.
\medskip

Let us begin by introducing the $p$-Wasserstein distance: Given two measures $\mu$ and $\nu$ on $\R^d$ with the same mass $\mu(\R^d)=\nu(\R^d)$, denote by $\mathrm{Cpl}(\mu,\nu)$ the set of measures $\pi$ on $\R^d\times\R^d$ such that for all $A,B\subset\R^d$ (with the convention that $A\subset\R^d$ will denote that $A$ is a Borel measurable subset of $\R^d$.):
\begin{align*}
    \pi(A\times\R^d)=\mu(A)&&\text{and}&&\pi(\R^d\times B)=\nu(B).
\end{align*}
The $p$-Wasserstein cost is then defined by
\begin{align*}
    \W_p^p(\mu,\nu):=\inf_{\pi\in\mathrm{Cpl}(\mu,\nu)}\int_{\R^d\times\R^d}|x-y|^p\pi(dx,dy).
\end{align*}

\subsection{Local transportation distances}\label{SS:1}

This subsection will deal with modifications of the Wasserstein distance to some domain $\Lambda\subset\R^d$. For two measures $\mu,\nu$ such that $\mu(\Lambda)=\nu(\Lambda)$, write
    \begin{align*}
        \W_\Lambda^p(\mu,\nu):=\W^p_p(\mu\llcorner \Lambda,\nu\llcorner \Lambda).
    \end{align*}
In the following, let $\mu,\nu$ be two locally finite measures and $\Lambda\subset\R^d$ a bounded open domain.
\medskip 

As mentioned in the introduction, the boundary transport's main distinction is that compared to the usual Wasserstein distance the measures need not have the same mass. This is accomplished by introducing an infinite mass reservoir at the boundary $\partial\Lambda$ of $\Lambda$:
\begin{align}\label{a1}
    \Wb^p_\Lambda(\mu,\nu)=\inf_{\pi\in\mathsf{Adm}_b(\mu,\nu)}\int_{\overline{\Lambda}\times\overline\Lambda}|x-y|^p\pi(dx,dy),
\end{align}
where the set of admissible couplings $\mathsf{Adm}_b(\mu,\nu)$ is the set of all measures $\pi$ with first marginal $\mu$ and second marginal $\nu$ in $\Lambda$. This means that
\begin{align*}
   \pi(A\times \overline\Lambda)=\mu(A) &&\text{and}&&\pi(\overline\Lambda\times A)=\nu(A) &&\text{for any }A\subset\Lambda.
\end{align*}
\medskip

The last local distance is the asymmetric transport $\Wa$, which has been introduced in \cite{DSS} to provide an alternative derivation for $\underline c(p,d)$. For technical reasons we however slightly modify the definition\footnote{In section \ref{S:4} we prove that our definition leads to the same constant. The only difference is that we again allow for measures with different mass, in the case of equal mass the definitions coincide.}:
\begin{align}\label{a2}
    \Wa^p_\Lambda(\mu,\nu)=\inf_{\lambda\le\mu}\W_{\overline\Lambda}^p(\lambda,\nu\llcorner\Lambda).
\end{align}
Here $\lambda\le\mu$ means that $\lambda(A)\le\mu(A)$ for all sets $A\subset\Lambda$. Note that $\lambda$ is allowed to charge mass to the boundary $\partial\Lambda$, as for the r.h.s. to be well-defined $\lambda(\overline\Lambda)=\nu(\Lambda)$. The rest of this section will be devoted to some properties of $\Wa^p$ as well as comparisons to the other distances. 
\begin{lemma}\label{L:opt}
    There exists a measure $\lambda\leq \mu$ on $\overline\Lambda$ and a coupling $\pi$ of $\lambda$ and $\nu\llcorner\Lambda$ such that
    \begin{align*}
        \Wa^p_\Lambda(\mu,\nu)=\int_{\overline{\Lambda}\times\overline{\Lambda}}|x-y|^pd\pi(x,y).
    \end{align*}
\end{lemma}
{\sc Proof.} Consider the set $\mathsf{Adm}_a(\mu,\nu)$ of measures $\pi$ on $\overline\Lambda\times\overline\Lambda$ which satisfy
\begin{align*}
    \pi(A\times \overline\Lambda)\le \mu(A) \text{ for all } A\subset\Lambda &&\text{and}&&\pi(\overline\Lambda\times B)=\nu(B)\text{ for all }B\subset{\Lambda}.
\end{align*}
This allows us to reformulate
\begin{align*}
    \Wa^p_\Lambda(\mu,\nu)=\inf_{\pi\in\mathsf{Adm}_a(\mu,\nu)}\int_{\overline\Lambda\times\overline\Lambda}|x-y|^p\pi(dx,dy).
\end{align*}
Since $\Lambda$ is a bounded set, $\pi(\overline\Lambda\times\overline\Lambda)=\nu(\Lambda)<\infty$, and $\mathsf{Adm}_a(\mu,\nu)$ is weakly compact. Combining with the lower semi-continuity of 
\begin{align*}
    \pi\longrightarrow \int_{\overline\Lambda\times\overline\Lambda}|x-y|^p\pi(dx,dy),
\end{align*}
existence of an optimizer follows. Defining $\lambda(\cdot)=\pi(\cdot\times\overline\Lambda)$ the claim follows.
\qed
\medskip

Due to the similarity of the definitions of $\Wa_\Lambda^p$ and $\Wb_\Lambda^p$ it is of course natural to compare them. Indeed, it turns out that the asymmetric transport distance is always smaller, as it does not require mass to be transported to the boundary but can rather just "appear" there. To make this precise, let $\pi\in\mathsf{Adm}_b(\mu,\nu)$. As mentioned in \cite{FG}, we can assume $\pi(\partial\Lambda,\partial\Lambda)=0$. Define from this the measure 
\begin{align*}
    \lambda(\cdot):=\tilde{\pi}(\cdot\times\overline{\Lambda}):=\big(\pi-\pi\llcorner(\Lambda\times\partial\Lambda)\big)(\cdot\times\overline\Lambda).
\end{align*}
By definition of $\pi\in\mathsf{Adm}_b(\mu,\nu)$, 
\begin{align*}
    \lambda(A)=\tilde{\pi}(A\times\overline\Lambda)\le\pi(A\times\overline\Lambda)=\mu(A)&&\text{for }A\subset\Lambda
\end{align*}
and $\tilde{\pi}$ is a coupling of $\lambda$ and $\nu\llcorner\Lambda$. Thus
\begin{align*}
    \W_{\overline\Lambda}^p(\lambda,\nu\llcorner\Lambda)\le \int_{\overline\Lambda\times\overline\Lambda}|x-y|^p\tilde{\pi}(dx,dy)\le\int_{\overline\Lambda\times\overline\Lambda}|x-y|^p\pi(dx,dy),
\end{align*}
which implies the desired estimate
\begin{align}\label{a3}
    \Wa_\Lambda^p(\mu,\nu)\le\Wb_\Lambda^p(\mu,\nu).
\end{align}
\subsection{Global transportation between random measures}\label{SS:2}
The drawback of the distances in the previous section is that they allow only for a bounded observation window or finite measures. The framework introduced in \cite{Hue, HS} allows to circumvent both of these restrictions in the case of random measures which satisfy a strong notion of stationarity. In the following we assume $(\Omega,\mathcal{F},\P) $ to be a probability space which supports a measurable flow $(\theta_z)_{z\in\R^d}$ such that $(\theta_z)_\#\P=\P$ for all $z\in\R^d$. For convenience, we will moreover assume $\Omega$ to be a compact standard Borel space. A random measure $\eta$ is a random variable on the space of locally finite measures on $\R^d$ equipped with the Borel $\sigma$-algebra and will be called equivariant if
\begin{align*}
    \eta^{\theta_z\omega}(A+z)=\eta^\omega(A) &&\text{for all }\omega\in\Omega, A\subset\R^d.
\end{align*}
Since the Lebesgue measure $\mathcal{L}$ is up to constants the unique translation invariant measure on $\R^d$, the intensity of an equivariant random measure is given by 
\begin{align*}
    \E[\eta(A)]=c\mathcal{L}(A) &&\text{for all }A\subset\R^d
\end{align*}
for some $c>0$. A prime example is the $c$-homogeneous Poisson point process $\eta$, which is defined by the fact the point statistic of each set $A$ is a Poisson random variable $\eta(A)\sim\mathsf{Poi}(c\mathcal{L}(A))$ and for a collection of disjoint sets $A_1,...,A_n$ the random variables $\eta(A_1),...,\eta(A_n)$ are independent.
\medskip

Let $\eta,\xi$ be two jointly equivariant measures\footnote{By this we mean that they are equivariant w.r.t. the same flow.} with equal intensity. By rescaling we may assume $c=1$. The class $\mathrm{Cpl}_\text{equiv}(\eta,\xi)$ consists of random measures $\Pi$ on $\R^d\times\R^d$, called equivariant couplings of $\eta$ and $\xi$, that satisfy for all $\omega\in\Omega$, $A,B\subset\R^d$:
\begin{align*}
    &&\Pi^\omega(A\times \R^d)=\eta^\omega(A),\Pi^\omega(\R^d\times B)=\xi^\omega(B)\\
    \text{and}&& \Pi^{\theta_z\omega}((A+z)\times(B+z))=\Pi^\omega(A\times B).
\end{align*}
With this definition at hand, we are in a position to define the mean transportation cost of a coupling $\Pi\in\mathrm{Cpl}_\text{equiv}(\eta,\xi)$:
\begin{align*}
    \mathsf{Cost}_p(\Pi):=\E\left[\int_{(0,1)^d\times\R^d}|x-y|^p\Pi(dx,dy)\right].
\end{align*}
Denoting $\Lambda_n:=(-\frac n2,\frac n2)^d$, $\Lambda_{n,z}:=\Lambda_n+z$ and using equivariance
we recover the identity
\begin{align*}
    \mathsf{Cost}_p(\Pi):=n^{-d}\E\left[\int_{\Lambda_{n,z}\times\R^d}|x-y|^p\Pi(dx,dy)\right]&&\text{ for all }n\in\N,z\in\R^d.
\end{align*}
The mean transport cost $\mathsf{Cost}_p(\eta, \xi)$ between $\eta$ and $\xi$ is the minimal cost over all equivariant couplings 
\begin{align}\label{D:cost}
    \mathsf{Cost}_p(\eta,\xi):=\inf_{\Pi\in\mathrm{Cpl}_\text{equiv}(\eta,\xi)}\mathsf{Cost}_p(\Pi).
\end{align}
A special case is given by $\eta$ a $1$-homogeneous Poisson point process and $\xi=\mathcal{L}$, where we define
\begin{align}\label{D:Global}
    c_\infty(p,d):=\mathsf{Cost}_p(\eta,\mathcal{L}).
\end{align}
We will prove a more general statement than Theorem \ref{T:2}, namely an extension to arbitrary equivariant random measures:
\begin{theorem}\label{T:1}
    For any two jointly equivariant random measures $\eta,\xi$ with equal intensity and any $d\in\N,p\ge 1$ both limits
    \begin{align*}
        \lim_{n\to\infty}n^{-d}\E[\Wa_{\Lambda_n}^p(\eta,\xi)] && \text{and} &&  \lim_{n\to\infty}n^{-d}\E[\Wb_{\Lambda_n}^p(\eta,\xi)]
    \end{align*}
    exist in $[0,\infty]$ and are equal to $\mathsf{Cost}_p(\eta,\xi)$.
\end{theorem}
This proves Theorem \ref{T:2} up to reintroducing the normalization of the Lebesgue measure, which we perform in Section \ref{S:4}.

\section{Proof of Theorem \ref{T:1}}
This section consists of proving the following inequalities:
\begin{proposition}\label{P:1}
For any jointly equivariant random measures $\eta,\xi$ with equal intensities
    \begin{align*}
        \limsup_{n\to\infty}n^{-d}\E[\Wb_{\Lambda_n}^p(\eta,\xi)]\le\mathsf{Cost}_p(\eta,\xi).
    \end{align*}
\end{proposition}

\begin{proposition}\label{P:2}
For any jointly equivariant random measures $\eta,\xi$ with equal intensity,
    \begin{align*}
        \mathsf{Cost}_p(\eta,\xi)\le \liminf_{n\to\infty}n^{-d}\E[\Wa^p_{\Lambda_n}(\eta,\xi)].
    \end{align*}  
\end{proposition}
Together with the estimate \eqref{a3} this proves Theorem \ref{T:1}.
\subsection{Construction of a competitor to the boundary transport}
The goal of this subsection is to prove Proposition \ref{P:1}. We assume $\mathsf{Cost}_p(\eta,\xi)<\infty$, as else the statement is obviously true. The idea of the proof is to use the optimizer of the $\mathsf{Cost}_p$-functional (see \eqref{D:cost} for the definition) to construct a competitor for the minimizing transport map of $\Wb_{\Lambda_n}(\eta,\xi)$. Let $\q$ be such an optimizer. 
\medskip

We begin by establishing that the cost of the mass being transported outside of $\Lambda_n$ is asymptotically negligible:
\begin{lemma}\label{L:1}
    As $n\to\infty$,
    \begin{align*}
        \E\left[\int_{\Lambda_n\times\Lambda_n^c}|x-y|^pd\q(x,y)\right],\E\left[\int_{\Lambda_n^c\times\Lambda_n}|x-y|^pd\q(x,y)\right]=o(n^d).
    \end{align*}
\end{lemma}
{\sc Proof.}
By symmetry it is of course enough to prove only one of the inequalities. We begin by dividing $\Lambda_n$ into a set close to the boundary and one at large distance:
\begin{align*}
    &\E\left[\int_{\Lambda_n\times\Lambda_n^c}|x-y|^pd\q(x,y)\right]\\
    &= \E\left[\int_{\Lambda_n\setminus\Lambda_{n-\sqrt{n}}\times\Lambda_n^c}|x-y|^pd\q(x,y)\right]+\E\left[\int_{\Lambda_{n-\sqrt{n}}\times\Lambda_n^c}|x-y|^pd\q(x,y)\right].
\end{align*}
Appealing to equivariance of $\q$, we estimate
\begin{align*}
    \E\left[\int_{\Lambda_n\setminus\Lambda_{n-\sqrt{n}}\times\Lambda_n^c}|x-y|^pd\q(x,y)\right]\le (n^d-(n-\sqrt{n})^d) \E\left[\int_{\Lambda_1\times\R^d}|x-y|^pd\q(x,y)\right],
\end{align*}
which is of order $n^{d-1/2}$. Again, by equivariance,
\begin{align*}
    \E\left[\int_{\Lambda_{n-\sqrt{n}}\times\Lambda_n^c}|x-y|^pd\q(x,y)\right]&\le\E\left[\int_{\Lambda_{n-\sqrt{n}}\times\R^d}1_{\{|x-y|>\sqrt{n}/2\}}|x-y|^pd\q(x,y)\right]\\
    &\le n^d\E\left[\int_{\Lambda_1\times\R^d}1_{\{|x-y|>\sqrt{n}/2\}}|x-y|^pd\q(x,y)\right].
\end{align*}
By monotone convergence, the expectation converges to  $0$ as $n\to\infty$, proving the statement.
\qed

\medskip

{\sc Proof of Proposition \ref{P:1}.}
In order to transform $\q$ into a competitor for $\Wb_{\Lambda_n}$, we define the projection $\mathsf{id}_{\Lambda_n}:\R^d\to\overline\Lambda_n$, such that $x\mapsto\arg\min_{y\in\overline{\Lambda}_n}|x-y|$. Then the measure $\pi=(\mathsf{id}_{\Lambda_n},\mathsf{id}_{\Lambda_n})_\#(\q\llcorner((\Lambda_n\times\R^d)\cup( \R^d\times\Lambda_n))\in\mathsf{Adm}_b(\eta,\xi)$, which implies
\begin{align*}
    \Wb_{\Lambda_n}^p(\eta,\xi)&\leq \int_{\overline{\Lambda}_n\times\Lambda_n\cup\Lambda_n\times\overline{\Lambda}_n}|x-y|^pd\pi(x,y)\\
    &\le \int_{\Lambda_n\times\R^d}|x-y|^pd\q(x,y)+\int_{\Lambda_n^c\times\Lambda_n}|x-y|^pd\q(x,y).
\end{align*}
Taking expectations and appealing to Lemma \ref{L:1}, we recover
\begin{align*}
    \E\left[\Wb^p_{\Lambda_n}(\eta,\xi)\right]&\le \mathsf{Cost}_p(\eta,\xi)|\Lambda_n|+\E\left[\int_{\Lambda_n^c\times\Lambda_n}|x-y|^pd\q(x,y)\right]\\
    &=\mathsf{Cost}_p(\eta,\xi)n^d+o(n^d).
\end{align*}

\qed

\subsection{Construction of an equivariant coupling}
Contrary to the subsection before, our goal is to construct an invariant coupling $\Pi$ of $\eta$ and $\xi$ on the whole space, such that $\mathsf{Cost}_p(\Pi)\le\liminf_{n\to\infty}n^{-d}\E[\Wa_{\Lambda_n}^p(\eta,\xi)]$. For this we follow the approach in \cite{Hue}.
\medskip

Throughout the whole section, we may assume 
\begin{align}\label{eq:liminfWa}
    \liminf_{n\to\infty}n^{-d}\E[\Wa^p_{\Lambda_n}(\eta,\xi)]<\infty, 
\end{align}
as else the statement would be trivially true. Moreover, we may assume that also 
\begin{align*}
    \limsup_{n\to\infty}n^{-d}\E[\Wa_{\Lambda_n}^p(\eta,\xi)]<\infty,
\end{align*}
as else we can pass to the subsequence realizing the $\liminf$.
\medskip

Let us begin by proving that the mass transported by an optimizer of the asymmetric transport restricted to $\Lambda_n$ through the boundary is of lower order compared to  $n^d$. 
\begin{lemma}\label{L:4}
    Let $\q_n^{\omega}$ be an optimizer of $\Wa_{\Lambda_n}(\eta^\omega,\xi^\omega)$. Then for any $0<\varepsilon<\frac p{p+1}$,
    \begin{align*}
        \E[\q_n(\partial\Lambda_n\times\Lambda_n)]=o(n^{d-\varepsilon}).
    \end{align*}
\end{lemma}
{\sc Proof.}
For $\alpha=\frac 1{p+1}$, we can estimate\footnote{By $A(n)\lesssim B(n)$ we mean that there exists a universal constant $C>0$ such that $A(n)\le CB(n)$ for all $n$. }
\begin{align*}
    &\E[\q_n(\partial\Lambda_n\times\Lambda_n)]\le\E[\q_n(\partial\Lambda_n\times(\Lambda_n\setminus \Lambda_{n-n^\alpha}))]+\E[\q_n(\partial\Lambda_n\times\Lambda_{n-n^\alpha})]\\
    &\lesssim\E[\q_n(\partial\Lambda_n\times(\Lambda_n\setminus \Lambda_{n-n^\alpha}))]+n^{-p\alpha}\E\left[\int_{\overline\Lambda_n\times \Lambda_{n-n^\alpha}}|x-y|^pd\q_n(x,y)\right]\\
    &\le(n^d-(n-n^\alpha)^d)+n^{-p\alpha}\E\left[\int_{\overline\Lambda_n\times \Lambda_{n}}|x-y|^pd\q_n(x,y)\right].
\end{align*}
Multiplying by $n^{-d+\varepsilon}$, we infer
\begin{align*}
    n^{-d+\varepsilon}\E[\q_n(\partial\Lambda_n\times\Lambda_n)]\lesssim n^{\alpha+\varepsilon-1}+n^{-p\alpha+\varepsilon}n^{-d}\E\left[\int_{\overline\Lambda_n\times\Lambda_n}|x-y|^pd\q_n(x,y)\right]
\end{align*}
By the choice of $\alpha$, both $\alpha+\varepsilon-1,\varepsilon-p\alpha<0$. The statement follows from 
\begin{align*}
    \limsup_{n\to\infty}n^{-d}\E\left[\int_{\Lambda_n\times\Lambda_n}|x-y|^pd\q_n(x,y)\right]<\infty.
\end{align*}
\qed

In the following we choose an optimizer $\q_n$ of $\Wa_{\Lambda_n}(\eta,\xi)$ with marginals $\lambda_n$ and $\xi\llcorner \overline\Lambda_n$, and define the measures
\begin{align}\label{a4}
    \q^\omega_{n,z}(\cdot\times\cdot):=\q_n^{\theta_{-z}\omega}((\cdot-z)\times(\cdot-z)).
\end{align}
From this we now construct a measure on the whole space:
\begin{align}\label{a5}
    \Pi_n^\omega(\cdot):=n^{-d}\sum_{z\in\Z^d}\q_{n,z}^\omega(\cdot).
 \end{align}
 In the following, we will adopt a slightly different viewpoint and view $\Pi_n$ as a measure $Q_n$ on $\R^d \times\R^d\times \Omega $. For some measurable $\Omega'\subset \Omega$, $A,B\subset \R^d$ we set:
 \begin{align*}
     Q_n(A\times B\times\Omega')=\int_{\Omega'}\Pi^\omega(A\times B)d\P(\omega).
 \end{align*}
 However, we will still denote
 \begin{align*}
     \Pi_n(A\times B\times\Omega)=\int_\Omega\Pi_n^\omega(A\times B)d\P(\omega)=\E[\Pi_n(A\times B)].
 \end{align*}
 
\begin{lemma}\label{L:6}
The sequence of measures $\{Q_n\}_{n\in\N}$ satisfies:
    \begin{enumerate}
   \item Each $\Pi_n$ is $\Z^d$-equivariant.
   \item $\{Q_n\}_n$ is vaguely tight as a sequence of measures on $\R^d\times\R^d\times\Omega$.
   \end{enumerate}
\end{lemma}

{\sc Proof.}
The first claim directly follows from the definition \eqref{a4} of $\q_{n,z}$. To prove the second claim, it is enough to show that $\E[\Pi_n(A)]$ is bounded uniformly in $n$ for each bounded $A\subset \R^d\times\R^d$. Indeed tightness in the vague topology follows from $\sup_n\Pi_n(K)<\infty$ for all $K\subset\R^d\times\R^d\times\Omega\text{ compact}$, see for example \cite{Kal}. However, since $\Omega$ was assumed to be compact, there is a compact $A\subset\R^d\times \R^d$ such that $K\subset A\times\Omega$. Moreover, for each bounded set $A$ there exists an odd $m\in\N$ such that $A\subset\R^d\times\Lambda_m$, we may restrict ourselves to these sets. Using the $\Z^d$-equivariance of $\Pi_n$, we can estimate
\begin{align*}
    \E[\Pi_n(\R^d\times\Lambda_m)]&\le\sum_{z\in\Z^d\cap\Lambda_m}\E[\Pi_n(\R^d\times\overline\Lambda_{1,z})]=m^d\E[\Pi_n(\R^d\times\overline\Lambda_1)]\\
    &\le m^dn^{-d}\sum_{z\in\Z^d\cap\overline\Lambda_n}\E[\q_{n,z}(\R^d\times\overline\Lambda_1)]\\
    &=m^dn^{-d}\sum_{z\in\Z^d\cap\overline\Lambda_n}\E[\q_n(\R^d\times\overline\Lambda_{1,z})]\\
    &\le m^d\big(\frac{n+2}n\big)^d\lesssim m^d,
\end{align*}
which completes the proof of the second claim.
\qed
\medskip

Due to tightness, the sequence $\{Q_n\}_n$ is vaguely compact and thus admits accumulation points. First restrict to a subsequence realising the $\liminf$ in \eqref{eq:liminfWa}. Then, let $Q_\infty$ be the limit of a further subsequence $\{Q_{n_k}\}_{k\in\N}$, which we relabel immediately. Let us define for each $z\in\Z^{2d}$,
\begin{align*}
    Q_{\infty,z}:=Q_\infty\llcorner \bigg(\big(\big[-\frac12,\frac12\big)^d\times \big[-\frac12,\frac12\big)^d\big)+z\bigg)\times\Omega
\end{align*}
so that $Q_\infty=\sum_{z\in\Z^{2d}}Q_{\infty,z}$. For each $z$, it is easy to see by vague convergence and a standard approximation, that $Q_{\infty,z}\big(\big([-\frac12,\frac12)^d\times [-\frac12,\frac12)^d\big)+z\big)\times\cdot\big)$ is absolutely continuous with respect to $\P$ and thus by the disintegration theorem, we may write
\begin{align*}
    Q_{\infty,z}=:\Pi_{\infty,z}^\omega d\P(\omega).
\end{align*}
This then defines $\Pi_\infty=\sum_{z\in \Z^{2d}}\Pi_{\infty,z}$ so that
\begin{align*}
     Q_{\infty}=\Pi_{\infty}^\omega d\P(\omega).
\end{align*}
\begin{proposition}\label{P:3}
  $\Pi_\infty$ is a $\Z^d$-equivariant coupling of $\eta$ and $\xi$.
\end{proposition}

{\sc Proof of Proposition \ref{P:3}.}
 $\Pi_\infty$ inherits $\Z^d$-equivariance from $\{\Pi_n\}$: For any continuous, compactly supported $f\in C_c(\R^d\times\R^d\times\Omega)$,
 \begin{align*}
     \int_{\R^d\times\R^d\times\Omega} f(x,y,\omega)&d\Pi_\infty(x,y)d\P(\omega)=\lim_{n\to\infty}\int_{\R^d\times\R^d\times\Omega} f(x,y,\omega)d\Pi_n^\omega(x,y)d\P(\omega)\\
     &=\lim_{n\to\infty}\int_{\R^d\times\R^d\times\Omega} f(x-z,y-z,\omega)d\Pi_n^{\theta_z\omega}(x,y)d\P(\omega)\\
     &=\int_{\R^d\times\R^d\times\Omega} f(x-z,y-z,\omega)d\Pi_\infty^{\theta_z\omega}(x,y)d\P(\omega).
 \end{align*}

The proof proceeds in the following way: We show that for all measurable sets $A$
\begin{align}\label{a6}
    \Pi_\infty^\omega(A\times\R^d)=\eta^\omega(A)&&\text{and}&&\Pi_\infty^\omega(\R^d\times A)=\xi^\omega(A)&&\P-\text{a.s..}
\end{align}
Since the Borel $\sigma$-algebra is countably generated, this suffices to deduce that $\Pi_\infty$ almost surely has marginals $\eta$ and $\xi$. Instead of \eqref{a6}, we will prove the following equivalent statement for the second marginal: For all $0\le\varphi\in C_c(\R^d\times\Omega)$,
\begin{align}\label{a6:2}
    \int_{\R^d\times\R^d\times\Omega}\varphi(y,\omega)d\Pi^\omega_\infty(x,y)d\P(\omega)=\int_{\R^d\times\Omega}\varphi(y,\omega)d\xi^\omega(y)d\P(\omega).
\end{align}
By definition \eqref{a5}, for each $n\in\mathbb{N}$,
\begin{align}\label{a9}
    \Pi_n^\omega(\cdot\times\R^d)\le\eta^\omega(\cdot)+r_n^\omega(\cdot) &&\text{ and }&&\Pi_n^\omega(\R^d\times \cdot)=\xi^\omega(\cdot),
\end{align}
$\P$-almost surely, where $r_n^\omega(\cdot)=n^{-d}\sum_{z\in\Z^d}\lambda^\omega_{n,z}(\cdot\cap\partial\Lambda_{n,z})$. In the following, let $0\le\varphi\in C_c(\R^d\times\Omega)$ such that $\supp(\varphi)\subset\Lambda_m\times\Omega$ and $(\psi_R)_R\subset C_c(\R^d)$ be an increasing approximation of $1_{\R^d}$ satisfying $0\le\psi_R\le 1$ and $\psi_R(x)=1$ for all $x\in\Lambda_R$. Appealing to  vague convergence
\begin{align*}
    &\left|\int_{\R^d\times\R^d\times\Omega}\varphi(y,\omega)d\Pi^\omega_\infty(x,y)d\P(\omega)-\int_{\R^d\times\Omega}\varphi(y,\omega)d\xi^\omega(y)d\P(\omega)\right|\\
    &=\lim_{R\to\infty}\left|\int_{\R^d\times\R^d\times\Omega}\varphi(y,\omega)\psi_R(x)dQ_\infty(x,y,\omega)-\int_{\R^d\times\Omega}\varphi(y,\omega)d\xi^\omega(y)d\P(\omega)\right|\\
    &=\lim_{R\to\infty}\lim_{n\to\infty}\left|\int_{\R^d\times\R^d\times\Omega}\varphi(y,\omega)\psi_R(x)dQ_n(x,y,\omega)-\int_{\R^d\times\Omega}\varphi(y,\omega)d\xi^\omega(y)d\P(\omega)\right|.
\end{align*}
Therefore, due to \eqref{a9}, 
\begin{align*}
    &\left|\int_{\R^d\times\R^d\times\Omega}\varphi(y,\omega)d\Pi^\omega_\infty(x,y)d\P(\omega)-\int_{\R^d\times\Omega}\varphi(y,\omega)d\xi^\omega(y)d\P(\omega)\right|\\
    &=\lim_{R\to\infty}\lim_{n\to\infty}\left|\int_{\R^d\times\R^d\times\Omega}\varphi(y,\omega)(\psi_R(x)-1)dQ_n(x,y,\omega)\right|\\
    &\le\|\varphi\|_\infty\lim_{R\to\infty}\lim_{n\to\infty}\E[\Pi_n(\Lambda_R^c,\Lambda_m)]\lesssim m^d\|\varphi\|_\infty\lim_{R\to\infty}\lim_{n\to\infty}\E[\Pi_n(\Lambda_{R-m}^c\times\overline\Lambda_1)]\\
    &\lesssim m^d\|\varphi\|_\infty\lim_{R\to\infty}R^{-p}\limsup_{n\to\infty}\E\left[\int_{\R^d\times\overline\Lambda_1}|x-y|^pd\Pi_n(x,y)\right]\\
    &=m^d\|\varphi\|_\infty\lim_{R\to\infty}R^{-p}\limsup_{n\to\infty}n^{-d}\E[\Wa_{\Lambda_n}^p(\eta,\xi)]=0,
\end{align*}
proving \eqref{a6:2}. To deduce \eqref{a6}, we will prove an intermediate result, namely
\begin{align}\label{a7}
    \int_{\R^d\times\R^d\times\Omega}\varphi(x,\omega)d\Pi_\infty^\omega(x,y)d\P(\omega)\le\int_{\R^d\times\Omega}\varphi(x,\omega)d\eta^\omega(x)d\P(\omega).
\end{align}
    This suffices, since \eqref{a7} implies
    \begin{align}\label{a8}
        \Pi_\infty^\omega(\cdot\times\R^d)\le\eta^\omega(\cdot)&&\P-\text{a. s.}.
    \end{align}
    
    Applying the mass transport principle \cite[Lemma 2.5]{Hue} to the function
    \begin{align*}
        f:\Z^d\times\Z^d\to \R, f(z_1,z_2)=\E[\Pi_\infty(\Lambda_{1,z_1}\times\Lambda_{1,z_2})]
    \end{align*}
    we recover the fact that
    \begin{align}\label{a10}
        \E[\Pi_\infty(\Lambda_1\times\R^d)]&=\sum_{z\in\Z^d}f(0,z)=\sum_{z\in\Z^d}f(z,0)\nonumber\\
        &=
        \E[\Pi_\infty(\R^d\times \Lambda_1)]=\E[\xi(\Lambda_1)]=1.
    \end{align}
    Note that we used $\Pi_\infty(\partial\Lambda_1\times\R^d)\le \eta(\partial\Lambda_1)=0$ and $\Pi_\infty(\R^d\times\partial\Lambda_1)= \xi(\partial\Lambda_1)=0$ almost surely. Using $\Z^d$-stationarity to reduce to $A\subset\Lambda_1$ and combining \eqref{a8} with \eqref{a10}, this proves \eqref{a6}. 
    \medskip
    
    To prove \eqref{a7} we estimate, using \eqref{a9} and fact that $0\le \psi_R\le 1$,
    \begin{align*}
    \int_{\R^d\times\R^d\times\Omega}&\varphi(x,\omega)d\Pi_\infty^\omega(x,y)d\P(\omega)=\lim_{R\to\infty}\int_{\R^d\times\R^d\times\Omega}\varphi(x,\omega)\psi_R(y)d\Pi_\infty^\omega(x,y)d\P(\omega)\\
    &=\lim_{R\to\infty}\lim_{n\to\infty}\int_{\R^d\times\R^d\times\Omega}\varphi(x,\omega)\psi_R(y)d\Pi_n^\omega(x,y)d\P(\omega)\\
    &\le \int_{\R^d\times\Omega}\varphi(x,\omega)\eta^\omega(x)d\P(\omega)+\lim_{n\to\infty}\int_{\R^d\times\Omega}\varphi(x,\omega)dr_n^\omega(x)d\P(\omega)\\
    &\le \int_{\R^d\times\Omega}\varphi(x,\omega)\eta^\omega(x)d\P(\omega)+\|\varphi\|_\infty\lim_{n\to\infty}\E[r_n(\Lambda_m)].
    \end{align*}
    By Lemma \ref{L:4}, $r_n(\Lambda_m)\to 0$ in $L^1$-norm, proving \eqref{a7}. 
\qed
\medskip

{\sc Proof of Proposition \ref{P:2}.}
For $\Pi_\infty$ to be an admissible competitor we have to stationarise the measure in the following way:
\begin{align*}
    \Pi^\omega(A\times B):=\int_{[0,1)^d}\Pi_\infty^{\theta_{-z}\omega}((A-z)\times(B-z))dz&&\text{ for all }\omega\in\Omega.
\end{align*}
The resulting measure $\Pi$ is now $\R^d$-stationary and inherits its marginals from $\Pi_\infty$, thus making it an admissible competitor. By lower semicontinuity and the Portmanteau theorem,
\begin{align*}
    \mathsf{Cost}_p(\Pi)&=\E\left[\int_{\Lambda_1\times\R^d}|x-y|^p d\Pi(x,y)\right]\\
    &=\int_\Omega\int_{[0,1)^d}\int_{(\Lambda_1-z)\times\R^d}|x-y|^pd\Pi^{\theta_{-z}\omega}_\infty(x,y)dzd\P(\omega)\\
    &=\E\left[\int_{\Lambda_1\times\R^d}|x-y|^p d\Pi_\infty(x,y)\right]\\
    &\le\liminf_{n\to\infty}\E\left[\int_{\Lambda_1\times\R^d} |x-y|^pd\Pi_n(x,y)\right]\\
    &=\liminf_{n\to\infty}n^{-d}\E\left[\int_{\overline\Lambda_n\times\R^d}|x-y|^pd\q_{n}(x,y)\right]\\
    &\le \liminf_{n\to\infty}n^{-d}\E[\Wa_{\Lambda_n}^p(\eta,\xi)].
\end{align*}
\qed

\section{Semidiscrete Matching}\label{S:4}
The goal of this section is to prove that our modifications, which allow us to consider measures with different mass, do not influence the value of the constants
\begin{align*}
    c_a(p,d):=\lim_{n\to\infty}n^{-d}\E[\Wa^p_{\Lambda_n}(\eta,\kappa_n\mathcal{L})]\\
     \underline c(p,d):=\lim_{n\to\infty}n^{-d}\E[\Wb^p_{\Lambda_n}(\eta,\kappa_n\mathcal{L})]
\end{align*}
when $\eta$ is a homogeneous Poisson point process. Here $\kappa_n:=\frac{\eta(\Lambda_n)}{n^d}$ is the normalization needed such that
\begin{align*}
    \eta(\Lambda_n)=\kappa_n\mathcal{L}(\Lambda_n).
\end{align*}
With this restriction to measures with the same mass the definition of boundary transport used in this work coincides with the one in \cite{GHT}.

\begin{proposition}\label{P:4}
    If $d\ge 3$, then
    \begin{align*}
        \lim_{n\to\infty}n^{-d}|\E[\Wb^p_{\Lambda_n}(\eta,\mathcal{L})]-\E[\Wb_{\Lambda_n}^p(\eta,\kappa_n\mathcal{L})]|=0
    \end{align*}
    and
    \begin{align*}
        \lim_{n\to\infty}n^{-d}|\E[\Wa^p_{\Lambda_n}(\eta,\mathcal{L})]-\E[\Wa_{\Lambda_n}^p(\eta,\kappa_n\mathcal{L})]|=0
    \end{align*}
\end{proposition}
{\sc Proof.}
We begin by proving the statement for the boundary transport and deduce the case of the asymmetric transport. By the Young inequality and the gluing lemma (see \cite{FG} for the boundary transport version), there exists $C>0$ such that for any $\varepsilon>0$,
\begin{align}\label{a11}
    \Wb^p_{\Lambda_n}(\eta,\mathcal{L})\le (1+\varepsilon)\Wb_{\Lambda_n}^p(\eta,\kappa_n\mathcal{L})+\frac{C}{\varepsilon^{p-1}}\Wb_{\Lambda_n}^p(\kappa_n\mathcal{L},\mathcal{L}).
\end{align}
Since the reverse inequality 
\begin{align}\label{a12}
    \Wb^p_{\Lambda_n}(\eta,\kappa_n\mathcal{L})\le (1+\varepsilon)\Wb_{\Lambda_n}^p(\eta,\mathcal{L})+\frac{C}{\varepsilon^{p-1}}\Wb_{\Lambda_n}^p(\kappa_n\mathcal{L},\mathcal{L})
\end{align}
holds by the same argument, it is enough to show 
\begin{align}\label{a13}
    \lim_{n\to\infty}n^{-d}\E[\Wb_{\Lambda_n}^p(\kappa_n\mathcal{L},\mathcal{L})]=0.
\end{align}
Indeed, assuming \eqref{a13} and combining \eqref{a11} and \eqref{a12}, we can estimate
\begin{align}\label{a15}
    \lim_{n\to\infty}&n^{-d}|\E[\Wb^p_{\Lambda_n}(\eta,\mathcal{L})]-\E[\Wb_{\Lambda_n}^p(\eta,\kappa_n\mathcal{L})]|\nonumber\\
    &\le \varepsilon\lim_{n\to\infty}n^{-d}\big(\E[\Wb_{\Lambda_n}^p(\eta,\mathcal{L})]+\E[\Wb_{\Lambda_n}^p(\eta,\kappa_n\mathcal{L})]+\frac{C}{\varepsilon^p}\E[\Wb_{\Lambda_n}^p(\kappa_n\mathcal{L},\mathcal{L})]\big)\nonumber\\
    &\le \varepsilon (2+\varepsilon)c_\infty(p,d)\xrightarrow{\varepsilon\to 0}0,
\end{align}
where the last inequality is due to another application of \eqref{a12}, Proposition \ref{P:1} and the fact that $c_\infty(p,d)<\infty$ for $d\ge 3$ can be found in \cite{HS}. Let us now prove \eqref{a13}: For this we work on the event $\kappa_n>1$, since the reverse $0<\kappa_n<1$ can be dealt with by rescaling. The case $\kappa_n=0$ can be treated very explicitly, since all the mass has to appear at the boundary and be distributed evenly. We define the map $T:\Lambda_n\to\partial\Lambda_n, x\mapsto\arg\min_{y\in\partial\Lambda_n}|x-y|$ so that for $\pi:=(\mathsf{id},T)_\#(\mathcal{L}\llcorner\Lambda_n)$ it holds that $\pi\in\mathsf{Adm}_b(\mathcal{L},0)$ and estimate
\begin{align*}
    \Wb^p_{\Lambda_n}(0,\mathcal{L})&=\Wb^p_{\Lambda_n}(\mathcal{L},0)\le\int_{\Lambda_n}|x-T(x)|^pdx\lesssim n^{d+p}.
\end{align*}
Therefore, 
\begin{align*}
    n^{-d}\E[\Wb^p_{\Lambda_n}(\kappa_n\mathcal{L},\mathcal{L})1_{\{\kappa_n=0\}}]\lesssim n^{p}\P(\eta(\Lambda_n)=0)=n^pe^{-n^d}=o(1),
\end{align*}
which is negligible. In the case $\kappa_n>1$, the idea is to identify a cube $Q_n$ such that $\kappa_n\mathcal{L}(Q_n)=\mathcal{L}(\Lambda_n)=n^d$, transport the mass from this cube to $\Lambda_n$ and the rest to the boundary. An easy calculation shows
\begin{align*}
    Q_n=\kappa_n^{-1/d}\Lambda_n
\end{align*}
 and we therefore define 
 \begin{align*}
     T:\Lambda_n\to\overline\Lambda_n, x\mapsto \begin{cases}
      \kappa_n^{1/d}x &\text{if }x\in Q_n\\
      \underset{y\in\partial\Lambda_n}{\arg\min}\:|x-y| & \text{otherwise}
    \end{cases};&&\text{and} &&\pi:=(\mathsf{id},T)_\#(\kappa_n\mathcal{L}\llcorner\Lambda_n).
 \end{align*}
 Then $\pi\in\mathsf{Adm}_b(\kappa_n\mathcal{L},\mathcal{L})$ and we can estimate
 \begin{align*}
     \Wb_{\Lambda_n}^p(\kappa_n\mathcal{L},\mathcal{L})&\le\int_{\overline{\Lambda}_n\times\overline{\Lambda}_n}|x-y|^p\pi(dx,dy)= \int_{\overline\Lambda_n}|x-T(x)|^pd\kappa_n\mathcal{L}(x)\\
     &=\kappa_n\int_{Q_n}|x-\kappa_n^{1/d}x|^pdx+\kappa_n\int_{\overline\Lambda_n\setminus Q_n}\mathrm{dist}(x,\partial\Lambda_n)^pdx\\
     &=\kappa_n\left(|1-\kappa_n^{1/d}|^p\int_{Q_n}|x|^pdx+\int_{\overline\Lambda_n\setminus Q_n}\mathrm{dist}(x,\partial\Lambda_n)^pdx\right)\\
     &\lesssim \kappa_n\left(|1-\kappa_n^{1/d}|^p\int_0^{\frac n2\kappa_n^{-1/d}}r^{d+p-1}dr+n^{d-1}\big(\frac n2-\frac n2\kappa_n^{-1/d}\big)^{p+1}\right)\\
     &\lesssim\kappa_n\left(|1-\kappa_n^{1/d}|^pn^{d+p}\kappa_n^{-1-p/d}+n^{d+p}(1-\kappa_n^{-1/d})^{p+1}\right)\\
     &=n^{d+p}(1-\kappa_n^{-1/d})^p\big(1+\kappa_n(1-\kappa_n^{-1/d})\big)\\
     &\le n^{d+p}(1-\kappa_n^{-1})^p(1+\kappa_n).
 \end{align*}
 In the case of $0<\kappa_n<1$, we can use the above calculation to estimate 
 \begin{align*}
    \Wb_{\Lambda_n}^p(\kappa_n\mathcal{L},\mathcal{L}) &=\kappa_n\Wb_{\Lambda_n}^p(\mathcal{L},\kappa_n^{-1}\mathcal{L})=\kappa_n\Wb_{\Lambda_n}^p(\kappa_n^{-1}\mathcal{L},\mathcal{L})\\
    &\lesssim\kappa_nn^{d+p}(1-\kappa_n^{1/d})^p\big(1+\kappa_n^{-1}(1-\kappa_n^{1/d})\big)\\
    &=n^{d+p}(1-\kappa_n^{1/d})^p(1+\kappa_n-\kappa_n^{1/d})\\
    &\le n^{d+p}(1-\kappa_n)^p.
 \end{align*}
 Let us now introduce the event $\Theta_n:=\{|\eta(\Lambda_n)-n^d|\le n^{\alpha d}\}$ for some $\frac12<\alpha<1$, which we choose later. On $\Theta_n$, thus the estimate
 \begin{align*}
     n^{-d}\Wb(\kappa_n\mathcal{L},\mathcal{L})\lesssim n^pn^{pd(\alpha-1)}
 \end{align*}
 holds, and choosing $\alpha$ such that $d\alpha-d<-1$ this term becomes of lower order. On the event $\Theta_n^c$ we use the rather crude estimate
 \begin{align*}
     n^{-d}\Wb_{\Lambda_n}^p(\kappa_n\mathcal{L},\mathcal{L})\le n^{-d}n^p\int_{\Lambda_n}|\kappa_n-1|dx=n^{p}|\kappa_n-1|.
 \end{align*}
 Combining this with the Bernstein inequality (see e.g. \cite{BLM}), for some $c>0$
 \begin{align*}
     n^{-d}\E[\Wb^p_{\Lambda_n}(\kappa_n\mathcal{L},\mathcal{L})]&\lesssim\P(\Theta_n)n^{p(1+\alpha d-d)}+n^p\E[|\kappa_n-1|1_{\Theta_n^c}]\\
     &\le n^{p(1+\alpha d-d)}+n^p\sqrt{\E[|\kappa_n-1|^2]\P(\Theta_n^c)}\\
     &\le n^{p(1+\alpha d-d)}+n^pn^{-d/2}\exp(-cn^{d(2\alpha-1)}).
 \end{align*}
 By our choice of $\alpha$, we recover 
 \begin{align*}
     \limsup_{n\to\infty}n^{-d}\E[\Wb^p_{\Lambda_n}(\kappa_n\mathcal{L},\mathcal{L})]=0.
 \end{align*}
 The main obstruction to arguing analogously in the asymmetric case is the fact that $\Wa_{\Lambda_n}^p$ does not define a distance and it is therefore unclear how to deduce the analogue of \eqref{a11} and \eqref{a12}. Here we present a different argument, which uses the comparison between $\Wa_{\Lambda_n}^p$ and $\Wb_{\Lambda_n}^p$: Let $\pi_1$ be an optimizer of $\Wa_{\Lambda_n}^p(\eta,\kappa_n\mathcal{L})$ and $\gamma$ be an optimizer of $\Wb_{\Lambda_n}^p(\kappa_n\mathcal{L},\mathcal{L})$ s.t. $\gamma\llcorner (\partial\Lambda_n\times\partial\Lambda_n)=0$. We set 
 \begin{align*}
     \pi_2:=\gamma\llcorner(\Lambda_n\times\overline\Lambda_n).
 \end{align*}
 By definition, the second marginal of $\pi_1$ and the first marginal of $\pi_2$ is $\kappa_n\mathcal{L}$, which allows us to apply the gluing lemma to "glue" together $\pi_1$ and $\pi_2$ along the second marginal to construct a measure $\pi$. We now define
 \begin{align*}
     \tilde\pi:=\mathrm{proj}_{1,3}\big(\pi\llcorner(\overline\Lambda_n\times\Lambda_n\times\Lambda_n)\big)+\gamma\llcorner (\partial\Lambda_n\times\Lambda_n).
 \end{align*}
 This way $\tilde\pi\in\mathsf{Adm}_a(\eta,\mathcal{L})$ and we thus may use it as a competitor in $\Wa_{\Lambda_n}^p(\eta,\mathcal{L})$. Using Young's inequality,
 \begin{align*}
     \Wa_{\Lambda_n}^p(\eta,\mathcal{L})&\le \int_{\overline\Lambda_n\times\overline\Lambda_n}|x-z|^pd \tilde\pi(x,z) \\
     &\le (1+\varepsilon)\int_{\overline\Lambda_n\times\overline\Lambda_n}|x-y|^pd\pi_1(x,y)+\frac{C}{\varepsilon^{p-1}}\int_{\overline\Lambda_n\times\overline\Lambda_n}|x-y|^pd\gamma(x,y)\\
     &=(1+\varepsilon)\Wa_{\Lambda_n}^p(\eta,\kappa_n\mathcal{L})+\frac{C}{\varepsilon^{p-1}}\Wb_{\Lambda_n}^p(\kappa_n\mathcal{L},\mathcal{L}).
 \end{align*}
 This estimate is sufficient to deduce the claim.
 Indeed, take the limit in $n$ and then send $\varepsilon$ to $0$ in the first step,  we have by \eqref{a3} in the third estimate, Proposition \ref{P:4} in the first equality and Theorem \ref{T:1} in the last two equalities
 \begin{align*}
     \liminf_{n\to\infty}n^{-d}\E[\Wa_{\Lambda_n}^p(\eta,\mathcal{L})]&\le \liminf_{n\to\infty}n^{-d}\E[\Wa_{\Lambda_n}^p(\eta,\kappa_n\mathcal{L})]\\
     &\le \limsup_{n\to\infty}n^{-d}\E[\Wa_{\Lambda_n}^p(\eta,\kappa_n\mathcal{L})]\\
     &\le \limsup_{n\to\infty}n^{-d}\E[\Wb_{\Lambda_n}^p(\eta,\kappa_n\mathcal{L})]\\
     &=\limsup_{n\to\infty}n^{-d}\E[\Wb_{\Lambda_n}^p(\eta,\mathcal{L})]= c_\infty(p,d)\\
     &=\liminf_{n\to\infty}n^{-d}\E[\Wa_{\Lambda_n}^p(\eta,\mathcal{L})].
 \end{align*}
 We therefore deduce that
 \begin{align*}
     \lim_{n\to\infty}n^{-d}\E[\Wa^p_{\Lambda_n}(\eta,\kappa_n\mathcal{L})]=\liminf_{n\to\infty}n^{-d}\E[\Wa^p_{\Lambda_n}(\eta,\mathcal{L})]=c_\infty(p,d).
 \end{align*}
 \qed
\medskip

Let us note that under suitable concentration assumptions, the proof generalizes to equivariant measures $\eta,\xi$. The calculation however becomes much more tedious and we leave it to the interested reader.

\section*{Acknowledgments}
The authors acknowledge funding by the DFG Priority Programme ``Random Geometric Systems'' (SPP 2265) - project P10 ``Optimal matching and balancing transport''. This research has been supported by the DFG through the Research Training Group ``Rigorous Analysis of Complex Random Systems'' (RTG 3027) funded by the Deutsche Forschungsgemeinschaft (DFG, German Research Foundation – Project Number 524444762) and under Germany's Excellence Strategy EXC 2044/2 –390685587, Mathematics Münster: Dynamics–Geometry–Structure.

\section*{AI Usage}
The mathematical arguments were developed by the authors without any assistance by AI. Large language models were used for proof and grammar checking, the text was however written solely by the authors.

\end{document}